\documentclass{article}
\usepackage{amsfonts}
\usepackage{amssymb}
\usepackage{latexsym}
\def\be{\begin{equation}}
\def\ee{\end{equation}}
\def\bea{\begin{eqnarray}}
\def\eea{\end{eqnarray}}
\def\qed{\hfill\mbox{$\Box$}\medskip}

\newcommand{\0}{^{\phantom{0}}}
\begin{document}
\begin{Large}
\centerline{\bf Every multiple of~$4$ except $212$, $364$, $420$, and $428$}
\centerline{\bf is the sum of seven cubes}
\end{Large}
\vspace*{2ex}
\centerline{Kent D.~Boklan and Noam D.~Elkies}
\vspace*{2ex}
\centerline{February, 2008; revised February, 2009}

\vspace*{5ex}

\begin{quote}
{\bf Abstract.}  It is conjectured that every integer $N>454$
is the sum of seven nonnegative cubes.  We prove the conjecture
when $N$\/ is a multiple of~$4$.
\end{quote}

{\large\bf 1 Introduction}

Waring famously asserted in his Meditationes Algebraic\ae\ of 1770:\footnote{
  Ellison [1971, p.10] reports that this statement appears on pages 203--204
  of the 1770 edition.
  In the English translation \cite{Waring} of the 1782 edition,
  this statement appears on page 336 as part~(9) of ``Theorem 9'';
  see also the discussion on page 379 of the same translation,
  and the first section of \cite{Ram07} for a fuller treatment of
  the history of this problem than we give here.

  \vspace*{2ex}

  }

{\it Omnis integer numerus vel est cubus, vel e duobus, tribus,
4, 5, 6, 7, 8, vel novem cubis compositus \ldots}

meaning that every positive integer is the sum of at most nine positive cubes
(equivalently, of exactly nine nonnegative cubes).
A proof of this was given by Wieferich in 1909 (with an error later patched by
Kempner [1912]).  Landau then proved that only a finite set
requires nine cubes, and Dickson [1939] 
identified this set with $\{23, 239\}$.
Linnik then established, ineffectively, that only a finite set requires
eight cubes, and seven suffice after some point \cite{Linnik}.
In 1984, McCurley gave an effective proof of Linnik's result,
demonstrating that every integer larger than $\exp(\exp(13.94))$
is the sum of seven positive cubes~\cite{McCurley}.
This was recently reduced to $\exp(524)$ \cite{Ram07}
using an analytic sieve argument.

It is believed that the exceptional set for Linnik's seven cubes theorem is
\bea
  & \{15,\; 22,\; 23,\; 50,\; 114,\; 167,\; 175,\; 186,\; 212,\;
   \qquad\qquad
\label{eq:exc7} \\
  & \qquad\qquad
   231,\; 238,\; 239,\; 303,\; 364,\; 420,\; 428,\; 454\}.
\nonumber
\eea
(Indeed one expects that every sufficiently large integer is the sum of
four positive cubes \cite{DHL}, but even such a statement with
four replaced by six is well beyond our ability to prove.)

In the other direction, it is shown in~\cite[Theorem~1 and Lemma~3]{BRZ} that
if $454 < N < 2.5 \cdot 10^{26}$ then $N$\/ is the sum of cubes of seven
nonnegative integers; and it is observed in \cite[p.60]{Ram07} that
the computation reported in \cite[433--434]{DHL}
raises the upper bound to $\exp(78.7) > 10^{34}$.
But the computation to raise this bound to $\exp(524)$,
and thus prove that (\ref{eq:exc7}) is the full exceptional set,
remains utterly infeasible.

We give a different kind of partial result, where $N$\/ is restricted
by a congruence condition but not by size:

{\bf Theorem.}
{\em Every multiple of $4$ except $212$, $364$, $420$, and $428$
can be written as the sum of seven nonnegative cubes.}

This is not the first such result, but the only earlier work
in this direction that we know of is the proof in \cite{BRZ} that
if $N \equiv 0$ or $\pm 1 \bmod 9$ and $N$\/ is an invertible
cubic residue mod~$37$ then $N$\/ is the sum of seven nonnegative cubes.
Bertault et al.\ note that $37$ could be replaced by various
larger primes congruent to~$1$ mod~$3$.  But the condition mod~$9$
is essential, and restricts $N \bmod 9$ to the three most common residues
for a sum of seven cubes.\footnote{
  The following table gives the distribution mod~$9$ of
  $N = n_1^3 + n_2^3 + \ldots + n_7^3$ among the $9^7$ possibilities of
  $(n_1,n_2,\ldots,n_7) \bmod 9$:
  $$
  \begin{array}{c||c|c|c|c|c}
  N \bmod 9 & 0 &  \pm 1  &  \pm 2  &  \pm 3  & \pm 4
  \\ \hline
  \rm{proportion}
      & 17.97\% & 16.32\% & 12.21\% &  7.68\% & 4.80\%
  \end{array}
  $$
  Note that of the seventeen exceptions listed in~(\ref{eq:exc7}),
  eleven are congruent to $\pm 4 \bmod 9$,
  and the remaining six to $-3 \bmod 9$.
  These rare residues $\pm 3$, $\pm 4$ are also the least easily
  accessible to current approaches to the seven-cube problem,
  including ours: it will be seen that we must work hardest to
  prove our theorem for $N$\/ congruent to $\pm 3$ or $\pm 4 \bmod 9$.
  Note that while our condition $4|N$\/ also puts $N$\/ in
  a probabilistically favored congruence class, the discrepancy
  is minuscule: a random sum of seven cubes is divisible by~$4$
  with probability only $25.39\%$.
  }
Compared with \cite{BRZ},
our new ingredient is the use of a quadratic form
\hbox{$Q = \sum_{i=1}^3 c_i\0 X_i^2$}
with $(c_1,c_2,c_3) \neq (1,1,1)$ which is nevertheless known
to represent all positive integers in certain arithmetic progressions.

The rest of the paper is organized as follows.
We review the basic identity~(\ref{eq:sum6}) for writing
suitable integers as sums of six integer cubes, then give a criterion
(Proposition~1) under which the cubes are positive.  The criterion
requires an auxiliary prime $p \equiv 2 \bmod 3$ in an interval
$(AN^{1/3},BN^{1/3})$ and a small enough positive integer~$x_0$
such that $p | N - x_0^3$.  We then choose $c_1,c_2,c_3$ and show
that a suitable $x_0$ exists provided $p$ satisfies a congruence condition
mod~$72$.  Finally we use the explicit bounds of~\cite{RamRum}
on the distribution of primes in arithmetic progressions,
plus a short further computation of prime chains, to prove
that such~$p$ exists if $N > N_0 = 10^{18}$.
This completes the proof because that $N_0$ is well below
the threshold of $2.5 \cdot 10^{26}$ of~\cite{BRZ}.

Since we formulate the proof so as to use only positive rather than
nonnegative cubes, we automatically get a representation of $N$\/
as a sum of seven positive cubes for all $N>10^{18}$ divisible by~$4$.
For $N \leq 10^{18}$, there are cases not listed in (\ref{eq:exc7})
for which $N$\/ is a sum of six or fewer nonnegative cubes
but not exactly seven; the largest of these is apparently $2408$
(as it happens a multiple of~$4$), whose only representations
as a sum of seven or fewer positive cubes are

\bea
2408 &\!\! = \!\!& 10^3 + 10^3 + 7^3 + 4^3 + 1^3
\nonumber \\
 &\!\! = \!\!& 12^3 + 8^3 + 5^3 + 3^3 + 2^3 + 2^3
\nonumber \\
 &\!\! = \!\!& 11^3 + 8^3 + 6^3 + 6^3 + 5^3 + 2^3
\nonumber \\
 &\!\! = \!\!& 10^3 + 10^3 + 6^3 + 4^3 + 4^3 + 4^3
\label{eq:2408}
\eea
and permutations of these four sums.

\vspace*{2ex}

{\large\bf 2 A six cube identity}

The use of the identities equivalent to the following Lemma
and its Corollary to study sums of seven cubes
goes back at least as far as \cite[\S2]{Linnik}.

{\bf Lemma 1.} {\em
Let
\be
Q = \sum_{i=1}^3 c_i\0 X_i^2,
\qquad
C = \sum_{i=1}^3 c_i^3.
\label{eq:C,Q}
\ee
Then
\be
\sum_{i=1}^3 \left( (c_i p + X_i)^3 + (c_i p - X_i)^3 \right)
= 2 C p^3 + 6 p Q.
\label{eq:sum6}
\ee
}

{\em Proof}\/: Apply the identity $(r+s)^3 + (r-s)^3 = 2r^3 + 6rs^2$
to each of the three terms in the sum.\qed

{\bf Corollary.} {\em
If $p$ and the $c_i$ are positive integers
and the $X_i$ are integers such that $|X_i| < c_i p$ for each~$i$,
then $2 C p^3 + 6 p Q$ is the sum of cubes of six positive integers.
}\qed

\vspace*{2ex}

{\large\bf 3 Strategy}

Given $N$, we choose integers $c_i>0$ and
a prime \hbox{$p \equiv 2 \bmod 3$} of size roughly $N^{1/3}$,
and let $x_0$ be a small nonnegative integer such that
$N - x_0^3 = 2 C p^3 + 6 p Q_0$ for some integer $Q_0 \geq 0$.
Such an $x_0$ exists if $N/p^3$ is large enough
because every integer is a cube mod~$6p$.
If, moreover, $x_0$ can be chosen such that
$Q_0$ is represented by the quadratic form~$Q$,
and $N/p^3$ is small enough so that all the terms
in the sum in the identity (\ref{eq:sum6}) are positive,
then we can use that identity to write $N$\/ as
the sum of cubes of seven positive integers.

Because no ternary quadratic form represents all $Q_0 \geq 0$,
we may need to put a further congruence condition on~$x_0$
modulo some $\beta$ relatively prime to~$6p$; we shall then choose
the least $x_0 > 0$ satisfying these congruences, so that
$x_0 \leq 6 \beta p$.  We find that this imposes the following
lower and upper bounds on $N/p^3$:

{\bf Proposition 1} {\em
Let $c_1,c_2,c_3$ be positive integers with $c_1 = \min {c_i}$.
Set $C = \sum_{i=1}^3 c_i^3$. For some $\beta \geq 1$ assume that
\be
2C + 216 \beta^3 < \frac{N}{p^3} < 2C + 6 c_1^3.
\label{eq:ineq}
\ee
If $x_0$ is a positive integer such that $x_0 \leq 6\beta p$ and
\be
N - x_0^3  =  2 C p^3 + 6 p Q_0,
\label{eq:Q0}
\ee
then $Q_0 > 0$; if further
\be
Q_0\0 = c_1\0 X_1^2 + c_2\0 X_2^2 + c_3\0 X_3^2
\label{eq:X}
\ee
for some integers $X_i$, then $N$\/ is the sum of cubes of
seven positive integers.
}

{\em Proof}\/: The lower bound on~$N/p^3$ assures that
$$
N - x_0^3 \geq N - (6\beta p)^3 > 2 C p^3,
$$
so $Q_0 > 0$.  Given a solution of (\ref{eq:X}),
we use the identity (\ref{eq:sum6}) to write $N-x_0^3$ as
the sum of cubes of six integers.  Since $x_0 \geq 0$,
it thus suffices to verify that $|X_i| < c_i p$.  Indeed we have
$$
c_i X_i^2 \leq Q_0 = \frac{N - x_0^3 - 2C p^3}{6p}
  \leq \frac{N - 2C p^3}{6p}
  < \frac{6c_1^3 p^3}{6p} = c_1^3 p^2 \leq c_i^3 p^2
$$
so $X_i^2 < c_i^2 p^2 = (c_i p_i)^2$.
Since $c_i p > 0$, we are done.\qed

The inequalities (\ref{eq:ineq}) require
\be
\beta < c_1 / \root 3 \of{36}.
\label{eq:beta}
\ee
Assuming this condition holds, (\ref{eq:ineq}) restricts $p$ to
an interval $(A N^{1/3}, B N^{1/3})$ for some constants
$A,B$ with $0 < A < B$.  We then use explicit bounds
for the distribution of primes in arithmetic progressions to find
$N_0 < 10^{26}$ such that a suitable $p$ exists for all $N \geq N_0$.

When $N$\/ is a multiple of~$4$ but not of~$8$,
we shall need to impose a condition on $p \bmod 8$;
when $N$\/ falls in one of the hardest residue classes $\pm 4 \bmod 9$,
we shall also impose a condition on~$p \bmod 9$.

\vspace*{2ex}

{\large\bf 4 Choices and analysis}

We choose $(c_1,c_2,c_3) = (4\beta,4\beta,6\beta)$ where
$\beta = 1$ or $\beta = 5$ depending on the residue of $N \bmod 9$
(as specified later in this section).
Then condition (\ref{eq:beta}) is satisfied, and
\be
c_1\0 X_1^2 + c_2\0 X_2^2 + c_3\0 X_3^2 = 2 \beta (2X_1^2 + 2X_2^2 + 3X_3^2).
\label{eq:betaQ}
\ee
We calculate $C = 344 \beta^3$, whence
$A = \beta^{-1}/\root3\of{1072}$ and
$B = \beta^{-1}/\root3\of{904}$,
with ratio $B/A = (134/113)^{1/3} > 1.0584$.

We choose $x_0$ so that
\be
x_0^3 \equiv N - 2 C p^3 = N - 688(\beta p)^3 \bmod 6 \beta p,
\label{eq:x_0}
\ee
which is possible because every integer is a cube mod~$6 \beta p$.
We select the least positive~$x_0$ satisfying this congruence;
thus $x_0 \leq 6 \beta p$.
As $N$\/ and $6 \beta p$ are even, so is $x_0$.
Since also $4|N$\/ and $2|C$, while $\beta$ and $p$ are odd,
it follows that $N - 2Cp^3$ is a multiple of $12 \beta p$.
That is,
\be
Q_1 := \frac{N - x_0^3 - 688 (\beta p)^3}{12 \beta p}
\label{eq:Q1}
\ee
is an integer.

This integer $Q_1$ is positive by our Proposition.
Set $Q_0 = 2 \beta Q_1$.
In view of~(\ref{eq:betaQ}), we need $Q_1$ to be represented
by the quadratic form $2 X_1^2 + 2 X_2^2 + 3 X_3^2$.
This quadratic form is unique in its genus,
so it represents all nonnegative integers
that are not excluded by congruence conditions.  In this case
this means all $Q_1$ that are neither congruent to $1 \bmod 8$
nor of the form $9^t (9m+6)$ for some nonnegative integers $m,t$.
(See \cite[(16), pages 44--45]{Dickson} for this characterization
of the integers represented by $Q_1$.)  It remains to choose $p$
so that these conditions are satisfied.

The condition for $Q_1$ mod~$8$ holds automatically if $N \equiv 0 \bmod 8$,
because then the numerator $N - x_0^3 - 688 (\beta p)^3$ in~(\ref{eq:Q1})
is a multiple of~$8$ while the denominator is not,
so the quotient $Q_1$ is even.  Assume then that $N  \equiv 4 \bmod 8$.
Since $688 \equiv 16 \bmod 32$ and $\beta p$ is odd, we have
$688 (\beta p)^3 \equiv 16 \bmod 32$ as well.
Since $x_0$ is even, $x_0^3$ is
a multiple of~$8$ {\em not}\/ congruent to $16 \bmod 32$.
Therefore
\be
\frac{N - x_0^3 - 688 (\beta p)^3}{4} \not\equiv \frac{N}{4} \bmod 8.
\label{eq:Qmod8}
\ee
We therefore choose $p$ so that
\be
3 \beta p \equiv \frac{N}{4} \bmod 8,
\label{eq:pmod8}
\ee
and this guarantees that $Q_1 \not\equiv 1 \bmod 8$,
so $Q_1$ passes the \hbox{mod-$8$} test for representability by
the quadratic form $2 X_1^2 + 2 X_2^2 + 3 X_3^2$.

Next we ensure that $Q_1 \neq (9m+6) 9^t$ by choosing
$\beta \in \{1,5\}$ and $p \bmod 9$ so that
$Q_1 \not \equiv 0, 6 \bmod 9$.
Since $p \equiv 2 \bmod 3$ implies $p^3 \equiv -1 \bmod 9$,
the choice of $\beta$ determines $688(\beta p)^3 \bmod 9$.
Thus, as $688 \equiv 4 \bmod 9$, we have
\be
N - 688(\beta p)^3 \equiv N + 4 \beta^3 \bmod 9.
\label{eq:mod9}
\ee
Now $x_0^3 \equiv 0$ or $\pm 1 \bmod 9$ for all integers $x_0$.
Therefore if
\be
N + 4 \beta^3 \not\equiv 0 {\rm\ or\ } \pm\! 1 \bmod 9
\label{eq:easy9}
\ee
then $N - x_0^3 - 688 (\beta p)^3$ cannot be a multiple of~$9$,
whence $Q_1$ is not a multiple of~$3$,
so {\em a fortiori} not congruent to~$0$ or $6 \bmod 9$.
We have $\beta^3 \equiv 1 \bmod 9$ for $\beta = 1$,
and $\beta^3 \equiv -1 \bmod 9$ for $\beta = 5$;
thus, unless $N \equiv \pm 4 \bmod 9$, we may choose
$\beta$ so that $N$\/ satisfies condition (\ref{eq:easy9})

In the remaining cases $N \equiv \pm 4 \bmod 9$,
we choose $\beta$ so that $3|x_0$
by requiring $\beta = 5$ for $N \equiv 4 \bmod 9$
and $\beta = 1$ for $N \equiv -4 \bmod 9$.
Then $27 | x_0^3$, and the numerator $N - x_0^3 - 688 (\beta p)^3$
in~(\ref{eq:Q1}) is a multiple of~$9$ that we can control mod~$27$
by choosing $p \bmod 9$.  Indeed if $\beta p = 3k \pm 1$ then
$(\beta p)^3 \equiv 9k \pm 1 \bmod 27$. Since $688 \equiv 13 \bmod 27$
this gives
\be
\frac{N - x_0^3 - 688 (\beta p)^3}{9}
\equiv \frac{N \pm 13}{9} - k \bmod 3,
\label{eq:mod3}
\ee
with the sign chosen so that $9 | N \pm 13$.  Dividing by $4 \beta p$,
we conclude that
\be
\frac{Q_1}{3} \equiv \beta \left( k - \frac{N \pm 13}{9} \right) \bmod 3.
\label{eq:Q1mod3}
\ee
We thus choose $k \equiv \beta + ((N \pm 13) / 9) \bmod 3$. That is,
\be
\beta p \equiv 3\beta + \frac{N \pm 13}{3} \pm 1 \bmod 9,
\label{eq:bpmod9}
\ee
with the sign in $\pm 1$ chosen so that $\pm 1 \equiv -\beta \bmod 3$.
Then $Q_1 / 3 \equiv 1 \bmod 3$,
so $Q_1$ passes the \hbox{mod-$9$} test for representability by
$2 X_1^2 + 2 X_2^2 + 3 X_3^2$.

\vspace*{2ex}

{\large\bf 5 Conclusion}

To finish the proof of our Theorem, we show:

{\bf Lemma 2.} {\em
For $\beta \in \{1,5\}$ let
$A = \beta^{-1}/\root3\of{1072}$ and $B = \beta^{-1}/\root3\of{904}$.
Set $N_0 = 10^{18}$.
Then whenever $N > N_0$
there exists a prime $p \in (A N^{1/3}, B N^{1/3})$
in each odd congruence class $l \bmod 72 = 8 \cdot 9$
with $l \equiv 2 \bmod 3$.
}

This will suffice because $10^{18}$
is smaller than the lower bound of $2.5 \cdot 10^{26}$
of~\cite{BRZ} on an integer~$N$\/ not listed in~(\ref{eq:exc7})
that is not the sum of seven nonnegative cubes.

{\em Proof}\/:
Taking $k=72$ in \cite[Theorem~1]{RamRum}, we find that
for every prime $p>10^{10}$ there exists a prime $p' > p$ such that
$
 p' \equiv p \bmod 72$ and $p' \leq
 ((1+\epsilon_{72}\0)/(1-\epsilon_{72}\0)) \, p
$.
Consulting \cite[\S5, p.419, Table~1]{RamRum}, we find
$\epsilon_{72}\0 < 0.013$, so $(1+\epsilon_{72}\0)/(1-\epsilon_{72}\0) < 1.027$
which is well below the gap ratio of $B/A = (134/113)^{1/3} > 1.0584$
that we need.
This establishes Lemma~2 for $N > 1072 (5\cdot 10^{10})^3$.

That bound is not small enough for our application because
it exceeds the threshold of~\cite{BRZ}
and even (by a factor of~$13.4$) the improved bound of~$10^{34}$
reported in~\cite[p.60]{Ram07}.  But it reduces the proof of Lemma~2
to a finite computation.  To make this computation manageable, we prove:

{\bf Sublemma.} {\em
In each congruence class $l \bmod 72$ coprime to~$72$
there exist an integer $M_l$
and primes $p_i$ $(0 \leq i \leq M_l)$ such that $p_0 < 19541$,
$p_{i-1} < p_i < 1.0584 p_{i-1}$ for each $i=1,2,\ldots,M_l$,
and $p_{M_l}\0 > 10^{10}$.
}

To derive Lemma~2 from this Sublemma, let $p$ be the smallest prime such that
$p \equiv l \bmod 72$ and $p > A N^{1/3}$.
Then $p > 19541$ because
$$
N > N_0 = 10^{18} > 1072(5 \cdot 19541)^3.
$$
If $p > 10^{10}$, use~\cite{RamRum}.
Else apply the Sublemma and find the maximal nonnegative $i < M_l$
such that $p > p_i$.  Then
$$
p \leq p_{i+1} < 1.0584\ p_i < 1.0584\ p < 1.0584\ A N^{1/3} < B N^{1/3},
$$
and we are done.

(Conversely, the existence of the gap ratio $17573/16493 < 1.0655$
between the primes $16493$ and $17573$ congruent to $5 \bmod 72$
shows that $N_0$ cannot be brought much below $10^{18}$ in Lemma~1.)

To find the primes $p_i$ required by the Sublemma,
we run the following algorithm for each~$l$\/
with some choice of positive $\delta < 0.0584$:
\begin{itemize}
\item
  Let $p_0$ be the largest prime $p<19541$ such that
  $p \equiv l \bmod 72$.  Set $i=0$.
\item While $p_i < 10^{10}$, let $p_{i+1}$ be the least prime
  such that $p_{i+1} \equiv p_i \bmod 72$ and
  $p_{i+1} > (1+\delta) p_i$, and increment~$i$ to $i+1$.
\end{itemize}
Each ``largest prime'' and ``least prime'' is found by simply
stepping down or up the congruence class until a prime is found.
The second step is repeated at most
$\log(10^{10})/\log(1+\delta) \doteq 23/\delta$ times.
Once $p_i$ exceeds $10^{10}$, we succeed if
$\max_{i < M_l}\0 p_{i+1}/p_i < 1.0584$;
otherwise we try again for a smaller~$\delta$:
decreasing $\delta$ lengthens the computation but skips fewer primes.
We find that $\delta = 0.01$ is small enough for the computation
to succeed for each~$l$.  Then $M_l < 1250$ in each case,
and the largest $p_{i+1}/p_i$ ratio occurs at $(l,i)=(5,1)$,
namely $21101 / 19949 < 1.0578 < 1.0584$.
This computation, programmed in \hbox{{\sc pari-gp}~\cite{gp}},
takes less than a minute to run on an office desktop machine
(G5 PowerMac with a 1.8 GHz processor),
and completes the proof of the Sublemma and thus of our Theorem.

{\em Remark}\/:
We could also have used \cite[Corollary 5.2.2]{RamRum} to prove
a result similar to Lemma~2 but with the weaker bound $N_0 = 10^{21}$,
which is still good enough because $10^{21} < 2.5 \cdot 10^{26}$.
But this would implicitly rely on the much more extensive calculations
by Ramar\'{e} and Rumely of the distribution in arithmetic progressions
of the primes up to $10^{10}$,
which involve the computation of hundreds of millions of primes,
whereas we needed fewer than $1250 \, \phi(72) = 3 \cdot 10^4$ primes
to prove Sublemma~$2$.

\vspace*{2ex}

{\large\bf Acknowledgements}

We thank the referee for a careful reading of our manuscript
and for several helpful suggestions.  The first-named author thanks
Prof.~R.~C.~Vaughan for several helpful suggestions.
The second-named author's work is partly supported by the
National Science Foundation under grant DMS-501029;
he thanks Scott Kominers for pointers to the
literature on ternary quadratic forms.

\begin{small}

KDB:
\textsc{Department of Computer Science, Queens College,
Flushing,  NY 11367, U.S.A.} (\textsf{boklan@boole.cs.qc.edu})

NDE:
\textsc{Department of Mathematics, Harvard University,
Cambridge, MA 02138, U.S.A.} (\textsf{elkies@math.harvard.edu})

\end{small}

\end{document}